\documentclass[12pt]{amsart}

\newtheorem{theorem}{Theorem}[section]
\newtheorem{lemma}[theorem]{Lemma}

\theoremstyle{remark}
\newtheorem{remark}[theorem]{Remark}
\numberwithin{equation}{section}

\let \la=\lambda

\let \o=\omega
\let \a=\alpha
\let \f=\varphi
\let \b=\beta

\begin{document}
\title
{On modular inequalities in variable $L^p$ spaces}

\author{Andrei K. Lerner}
\address{Department of Mathematics,
Bar-Ilan University, 52900 Ramat Gan, Israel}
\email{aklerner@netvision.net.il}

\begin{abstract}
We show that the Hardy-Littlewood maximal operator and a class of
Calder\'on-Zygmund singular integrals satisfy the strong type
modular inequality in variable $L^p$ spaces if and only if the
variable exponent $p(x)\sim const$.
\end{abstract}

\keywords{Maximal function, singular integral, $A_{\infty}$
weight, variable~$L^p$}

\subjclass[2000]{Primary 42B20, 42B25}
\maketitle

\section{Introduction}
Let $p:{\mathbb R}^n\to [1,\infty)$ be a measurable function.
Denote by $L^{p(\cdot)}({\mathbb R}^n)$ the Banach space of
measurable functions $f$ on ${\mathbb R}^n$ such that for some
$\la>0$,
$$\int_{{\mathbb R}^n}|f(x)/\la|^{p(x)}dx<\infty,$$
with norm
$$\|f\|_{L^{p(\cdot)}}=\inf\left\{\la>0:\int_{{\mathbb
R}^n}|f(x)/\la|^{p(x)}dx\le 1\right\}.$$

The spaces $L^{p(\cdot)}({\mathbb R}^n)$ are a special case of
Musielak-Orlicz spaces (cf.~\cite{Mu}). The behavior of some
classical operators in harmonic analysis on $L^{p(\cdot)}({\mathbb
R}^n)$ is intensively investigated during several last years.
Among numerous papers appeared in this area, let us mention only
those of specific interest to us, to be precise those where
different aspects concerning the boundedness on
$L^{p(\cdot)}({\mathbb R}^n)$ of the Hardy-Littlewood maximal
operator \cite{CUFN,Di1,Di2,Ne1,Ne2,PR} and the Calder\'on-Zygmund
operators \cite{DR,KL} were studied.

We recall that the Hardy-Littlewood maximal operator is defined
for any $f\in L^1_{loc}({\mathbb R}^n)$ by
$$Mf(x)=\sup_{B\ni x}\frac{1}{|B|}\int_B|f(y)|dy,$$
where the supremum is taken over all balls $B$ containing $x$.

Let $p_-=\operatornamewithlimits{ess\, inf}_{x\in\mathbb{R}^n}
p(x)>1$ and $p_+=\operatornamewithlimits{ess\,
sup}_{x\in\mathbb{R}^n} p(x)<\infty.$ It has been proved by
Diening \cite{Di1} that if $p$ satisfies the following uniform
continuity condition
\begin{equation}\label{unif}
|p(x)-p(y)|\le\frac{c}{-\log|x-y|},\quad |x-y|<1/2,
\end{equation}
and if $p$ is a constant outside some large ball, then
\begin{equation}\label{h-l}
\|Mf\|_{L^{p(\cdot)}}\le c\|f\|_{L^{p(\cdot)}}
\end{equation}
for all $f\in L^{p(\cdot)}({\mathbb R}^n)$. After that the second
condition on $p$ has been improved independently in several
directions by Cruz-Uribe, Fiorenza, and Neugebauer \cite{CUFN} and
Nekvinda \cite{Ne1}. For example, it is shown in \cite{CUFN} that
if $p$ satisfies (\ref{unif}) and
$$
|p(x)-p(y)|\le \frac{c}{\log(e+|x|)},\quad |y|\ge |x|,
$$
then (\ref{h-l}) holds.

Diening and R\r{u}\v zi\v cka \cite{DR} (see also \cite[Theorem
2.7]{KL} and \cite[Section~8]{Di2}) have proved that if $p_->1$
and $p_+<~\infty$, then a large class of Calder\'on-Zygmund
operators is bounded on $L^{p(\cdot)}({\mathbb R}^n)$ provided the
Hardy-Littlewood maximal operator is bounded on
$L^{p(\cdot)}({\mathbb R}^n)$ and on $L^{(p(\cdot)/s)'}$ for some
$0<s<1$, where $p'(x)=p(x)/(p(x)-1)$.

A natural question arises about conditions on $p$ implying the
strong type inequality
\begin{equation}\label{gen}
\int_{{\mathbb R}^n}|Rf(x)|^{p(x)}dx\le c\int_{{\mathbb
R}^n}|f(x)|^{p(x)}dx
\end{equation}
(so-called modular inequality in terminology of Musielak
\cite{Mu}), where $R$ is any of the above-mentioned classical
operators. Note that in \cite{CUFN} the weak type modular
inequality
$$|\{x\in {\mathbb R}^n:Mf(x)>\a\}|\le c\int_{{\mathbb
R}^n}|f(x)/\a|^{p(x)}dx\quad(\a>0)$$ is proved under extremely
weak assumptions on $p$. It is easy to see that~(\ref{gen}) yields
the norm inequality
$$\|Rf\|_{L^{p(\cdot)}}\le c\|f\|_{L^{p(\cdot)}},$$
and therefore one should expect that the class of functions $p$,
for which~(\ref{gen}) holds, must be smaller than the
corresponding class implying~(\ref{h-l}). Nevertheless, our main
result is somewhat surprising, since it says that this class is
trivial. More precisely, we have the following.

\begin{theorem}\label{hl} {\it Let $p_->1$ and $p_+<\infty$. Then
the inequality
\begin{equation}\label{strong}
\int_{{\mathbb R}^n}(Mf(x))^{p(x)}dx\le c\int_{{\mathbb
R}^n}|f(x)|^{p(x)}dx
\end{equation}
holds for any $f\in L^{p(\cdot)}({\mathbb R}^n)$ if and only if
$p(x)\sim const$. }
\end{theorem}

It is noteworthy that analogous questions on singular integrals
are very similar to those when the boundedness on weighted
$L^p_{\o}$ implies $\o\in A_p$ (cf. \cite[p. 210]{Stein}). We
shall deal with a singular integral operator $Tf=f*K$, with kernel
$K$ satisfying the standard conditions
$$\|\widehat K\|_{\infty}\le c, \quad
|K(x)|\le c/|x|^n,$$
$$|K(x)-K(x-y)|\le c|y|/|x|^{n+1}\,\,\mbox{for}\,|y|<|x|/2,$$
and an additional nondegeneracy condition
$$|K(tu_0)|\ge c'/|tu_0|^n$$
for some unit vector $u_0$ and any $t\in {\mathbb R}$. Observe
that this class of operators contains, for instance, any one of
the Riesz transforms.

\begin{theorem}\label{sing} {\it Let $p_->1$ and $p_+<\infty$. Then
the inequality
\begin{equation}\label{strong1}
\int_{{\mathbb R}^n}|Tf(x)|^{p(x)}dx\le c\int_{{\mathbb
R}^n}|f(x)|^{p(x)}dx
\end{equation}
holds for any $f\in L^{p(\cdot)}({\mathbb R}^n)$ if and only if
$p(x)\sim const$. }
\end{theorem}

\section{Proofs}
By a weight we mean any non-negative locally integrable function
on ${\mathbb R}^n$. Given a ball $B$ and $f\in L^1_{loc}({\mathbb
R}^n)$, let $f_B=|B|^{-1}\int_Bf$. For measurable $f$ and $g$ the
notation $f\sim g$ means $f(x)=g(x)$ a.e.

We say that a weight $\o$ satisfies $A_{\infty}$ Muckenhoupt's
condition if for any $\a$, $0<\a<1$, there exists a $\b$,
$0<\b<1$, such that $|E|\ge \a|B|$ implies $\int_E\o dx\ge
\b\int_B\o dx$ for all balls $B$ and all subsets $E\subset B$.
There are many equivalent characterizations of $A_{\infty}$ (see,
e.g., \cite[Ch. 5]{Stein}). In particular, $\o\in A_{\infty}$ if
and only if (see \cite[p. 405]{GCRF} or \cite{Hr})
\begin{equation}\label{A}
\left(\frac{1}{|B|}\int_B\omega dx\right)
\exp\left(\frac{1}{|B|}\int_B\log(1/\omega)dx\right)\le A.
\end{equation}

We say that a family of weights $\{\o_{\a}\}_{\a\in {\mathcal A}}$
satisfies $A_{\infty}$ condition uniformly in $\a$ if $\o_{\a}\in
A_{\infty}$ for any $\a\in {\mathcal A}$ with corresponding
$A_{\infty}$ constants independing of $\a$.

\begin{lemma}\label{main}  {\it Let $p$ be a non-negative measurable function on
${\mathbb R}^n$. The family of weights $\{t^{p(x)}\}_{t>0}$
satisfies $A_{\infty}$ condition uniformly in $t$ if and only if
$p(x)\sim const$. }
\end{lemma}

\begin{proof}
When $p(x)\sim const$ the statement of the lemma is trivial. Thus,
we assume that $t^{p(x)}\in A_{\infty}$ uniformly in $t$. Applying
(\ref{A}) to $\o_t(x)=t^{p(x)}$ yields
\begin{equation}\label{res}
\frac{1}{|B|}\int_Bt^{p(x)-p_B}dx\le A
\end{equation}
for any ball $B$ and all $t>0$. Now, if $|\{x\in B:p(x)>p_B\}|>0$,
we get a contradiction by letting $t\to \infty$ in (\ref{res}).
Analogously, if $|\{x\in B:p(x)<p_B\}|>0$, we get a contradiction
by letting $t\to 0$ in (\ref{res}). Therefore, $p(x)=p_B$ for a.e.
$x\in B$ and for all balls $B$. Hence, the limit
$p_{\infty}=\lim_{|B|\to \infty}p_B,$ where it is taken over all
balls $B$ in ${\mathbb R}^n$ as the measure $|B|$ tends to
infinity, exists, and $p(x)=p_{\infty}$ for a.e. $x\in {\mathbb
R}^n$.
\end{proof}

We are now in a position to prove Theorems \ref{hl} and
\ref{sing}. Since for $p(x)\sim~const$ both theorems represent
known classical results, we need to prove only the converse
directions.

\begin{proof}[Proof of Theorem \ref{hl}.] It follows from
(\ref{strong}) that for any ball $B$ and any $f\in
L^1_{loc}({\mathbb R}^n)$,
\begin{equation}\label{inter}
\int_B(|f|_B)^{p(x)}dx\le c\int_B|f(x)|^{p(x)}dx.
\end{equation}
Let $E\subset B$ be an arbitrary measurable subset with $|E|\ge
\a|B|$, $0<\a<1$. Taking in (\ref{inter}) $f=t\chi_E$, $t>0$, we
get
$$\a^{p_+}\int_Bt^{p(x)}dx\le c\int_Et^{p(x)}dx.$$
Therefore, the family of weights $\{t^{p(x)}\}_{t>0}$ satisfies
$A_{\infty}$ condition uniformly in $t$. Now we invoke Lemma
\ref{main} to complete the proof.
\end{proof}

\begin{proof}[Proof of Theorem \ref{sing}.] We use the following
property of singular integrals (see \cite[Ch. 5, 4.6]{Stein}): for
any ball $B$ there exists a ball $B'$ such that $f_B\le
c|(Tf)\chi_{B'}|$ for any non-negative $f\in C_0^{\infty}$
supported in $B$ and $f_{B'}\le c|(Tf)\chi_{B}|$ for any
non-negative $f\in C_0^{\infty}$ supported in $B'$. It follows
from this and from (\ref{strong1}) that for such $f$,
\begin{equation}\label{inter1}
\int_{B'}(f_B)^{p(x)}dx\le c'\int_B(f(x))^{p(x)}dx
\end{equation}
and
\begin{equation}\label{inter2}
\int_B(f_{B'})^{p(x)}dx\le c'\int_{B'}(f(x))^{p(x)}dx.
\end{equation}
A simple limiting argument extends these estimates for any $f\ge
0$. Taking in (\ref{inter1}) $f=t\chi_E$, $t>0$, where $E\subset
B$ with $|E|\ge \a|B|$, we get
$$\a^{p_+}\int_{B'}t^{p(x)}dx\le c'\int_Et^{p(x)}dx.$$
However (\ref{inter2}), with $f=t\chi_{B'}$, yields
$$\int_Bt^{p(x)}dx\le c'\int_{B'}t^{p(x)}dx,$$
and therefore,
$$\a^{p_+}\int_Bt^{p(x)}dx\le c'^2\int_Et^{p(x)}dx.$$
This gives the desired result exactly as in the previous proof.
\end{proof}

\section{Concluding remarks}
\begin{remark} We recall that a weight $\o$ is doubling if
there exists a constant $c>0$ such that $\int_{2B}\o dx\le c
\int_{B}\o dx$ for any ball $B\subset {\mathbb R}^n$. It is well
known that any $A_{\infty}$ weight is doubling but the converse is
not true. In Lemma \ref{main}, the $A_{\infty}$ condition, in
general, can not be replaced by a wider doubling condition.
Indeed, one can construct on the real line disjoint sets $E_1$ and
$E_2$ of positive measure whose union is ${\mathbb R}^1$, while
$\chi_{E_1}$ and $\chi_{E_2}$ are doubling measures (see \cite[Ch.
1, 8.8]{Stein}). Let now $p(x)=c_1\chi_{E_1\times {\mathbb
R}^{n-1}}+c_2\chi_{E_2\times {\mathbb R}^{n-1}}$, where
$c_1\not=c_2$. Then $p(x)\not\sim const$, while it is easy to
check that $\int_{2B}t^{p(x)}dx\le c\int_{B}t^{p(x)}dx$ for any
ball $B\subset {\mathbb R}^n$ and all $t>0$.

However, assuming additionally that $p$ is continuous, one can
show that the family $\{t^{p(x)}\}_{t>0}$ is doubling uniformly in
$t$ if and only if $p(x)=const$.
\end{remark}

\begin{remark}
Musielak-Orlicz spaces (cf.~\cite{Mu}) consist of all measurable
$f$ such that for some $\la>0$,
$$\int_{{\mathbb R}^n}\f(x,|f(x)/\la|)dx<\infty,$$
where $\varphi:{\mathbb R}^n\times {\mathbb R}_+\to{\mathbb R}_+$
satisfies specific conditions. These spaces contain, i.e.,
weighted Lebesgue spaces $L^p_{\o}$ (when $\f(x,\xi)=\xi^p\o(x)$)
and Orlicz spaces (when $\f(x,\xi)$ is constant in the first
variable).

Theorems \ref{hl} and \ref{sing} show that in the case
$\f(x,\xi)={\xi}^{p(x)}$ the corresponding modular inequality for
$M$ or $T$ holds iff $\f$ is constant in the first variable. It is
easy to see that in general an analogous result does not hold. For
example, one can take $\f(x,\xi)=\xi^p\o(x)$ with $\o$ satisfying
the $A_p$ Muckenhoupt condition.

On the other hand, it is known (see, e.g., \cite{Ga}) that in the
context of Orlicz spaces the modular inequality for $M$ is
equivalent to the norm inequality. Theorem \ref{hl} shows that
this is not the case in the context of Musielak-Orlicz spaces.
\end{remark}

\vskip3mm
\noindent {\bf Acknowledgement.} The author thanks
Alexei Karlovich for helpful discussions.

\end{document}